\RequirePackage{fix-cm}
\documentclass[smallextended]{svjour3}       
\smartqed  

\usepackage{amsmath,graphicx,amssymb,amsfonts}
\newtheorem{thm}{Theorem}[section]
\newtheorem{cor}[thm]{Corollary}

\newtheorem{lem}[thm]{Lemma}

\newtheorem{defn}[thm]{Definition}
\newtheorem{exam}[thm]{Example}

\numberwithin{equation}{section}

\begin{document}

\title{Maximum nullity of Cayley graph
}


\author{ Ebrahim Vatandoost \and Yasser Golkhandy Pour}


\institute{E. Vatandoost \at
              Department of mathematics, Faculty of sciences, Imam Khomeini International University, Qazvin \\
              \email{vatandoost@si.ikiu.ac.ir}           
           \and
           Y. Golkhandy Pour \at
              Department of mathematics, Faculty of sciences, Imam Khomeini International University, Qazvin\\
               \email{y.golkhandypour@ikiu.ac.ir}
}

\date{Received: date / Accepted: date}

\maketitle

\begin{abstract}
One of the most interesting problems on maximum nullity (minimum rank) is to characterize $M(\mathcal{G})$ ($mr(\mathcal{G})$) for a graph $\mathcal{G}$. In this regard, many researchers have been trying to find an upper or lower bound for the maximum nullity. For more results on this topic, see \cite{4}, \cite{2}, \cite{10} and \cite{1}.
\newline
In this paper, by using a result of Babai \cite{Babai}, which presents the spectrum of a Cayley graph in terms of irreducible characters of the underlying group, and using representation and character of groups, we give a lower bound for the maximum nullity of Cayley graph, $X_S(G)$, where $G=\langle a\rangle$ is a cyclic group, or $G=G_1\times \cdots\times G_t$ such that $G_1=\langle a\rangle$ is a cyclic group and $G_i$ is an arbitrary finite group, for some $2\leq i\leq t$, with determine the spectrum of Cayley graphs.
\keywords{Cayley graphs\and Spectra of graphs\and Maximum nullity}
 \subclass{MSC 05C50.}
\end{abstract}

\section{Introduction}
For a positive integer $n$, let $S_n(\mathbb{R})$ be the set of all symmetric matrices of order $n$ over the real number. Suppose that $A\in S_n(\mathbb{R})$. Then the graph of $A$ which is denoted by $\mathcal{H}(A)$ is a graph with the vertex set $\{u_1, \ldots, u_n\}$ and the edge set $\{u_i\sim u_j : a_{ij}\neq 0, 0\leq i<j\leq n\}$. It should be noted that the diagonal of $A$ has no role in the determining of $\mathcal{H}(A)$.

The {\it set of symmetric matrices} of graph $\mathcal{G}$ is the set $S(\mathcal{G})=\{A\in S_n(\mathbb{R}) : \mathcal{H}(A)=\mathcal{G}\}$. The {\it minimum rank} of a graph $\mathcal{G}$ of order $n$ is defined to be the minimum cardinality between the rank of symmetric matrices in $S(\mathcal{G})$ and denoted by $mr(\mathcal{G})$.
Similarly, the {\it maximum nullity} of $\mathcal{G}$ is defined to be the maximum cardinality between the nullity of symmetric matrices in $S(\mathcal{G})$; and is denoted by $M(\mathcal{G})$. Clearly, $mr(\mathcal{G})+M(\mathcal{G})=n$.

One of the most interesting problems on minimum rank is to characterize $mr(\mathcal{G})$ for graphs. In this regard, many researchers have been trying to find an upper or lower bound for the minimum rank. For more results on this topic, see \cite{4}, \cite{2}, \cite{10} and \cite{1}.

The adjacency matrix of a graph $\mathcal{G}$ is the matrix $A_{\mathcal{G}}$ whose the entry $a_{ij}=1$ if and only if vertices $u_i$ and $u_j$ are adjacent, and $a_{ij}=0$ otherwise. The eigenvalues of $\mathcal{G}$ are the eigenvalues of $A_{\mathcal{G}}$, and the spectrum of $\mathcal{G}$ is the collection of its eigenvalues together with multiplicities. If $\lambda_1\ldots,\lambda_t$ are distinct eigenvalues of a graph $\mathcal{G}$ with respective multiplicity $n_1,\ldots n_k$, then we denote the spectrum of $\mathcal{G}$ by
\begin{equation}
 spec(\mathcal{G})=\bigg[\lambda_1^{n_1}, \ldots, \lambda_t^{n_t}\bigg].
\end{equation}
Let $G$ be a group, and let $S$ be a subset of $G$ that is closed under taking inverse and does not contain the identity, $e$. Then the {\it Cayley graph}, $X_S(G)$, is the graph with vertex set $G$ and edge set
\begin{equation}
  E=\{g_1\sim g_2\;:\; g_1g_2^{-1}\in S\}.
\end{equation}
Since $S$ is inverse-closed and does not contain the identity, it is a simple fact that $X(G,S)$ is undirected and has no loop.

In $1979$, Babai \cite{Babai} presented the spectrum of a Cayley graph in terms of irreducible characters of the underlying group $G$. The following important theorem was the result of this paper.
\begin{thm}\label{Baba} \cite{Babai}Let $G$ be a finite group of order $n$ whose irreducible characters (over
$\mathbb{C}$) are $\chi_{1},\ldots, \chi_{h}$ with respective
degree $n_{1},\ldots,n_{h}$. Then the spectrum of the Cayley graph
$X_S(G)$ can be arranged as $\Lambda = \{\lambda_{ijk} \;:\;
i=1 ,\ldots, h ; j,k =1 ,\ldots, n_i\}$ such that $\lambda_{ij1}=
\ldots=\lambda_{ijn_i}$ (this common value will be denoted by
$\lambda_{ij}$), and
\begin{equation}
\lambda_{i1}^{t}+\ldots+\lambda_{in_i}^{t}=\sum_{s_1,\ldots,s_t\in
S} \chi_i(\prod_{l=1}^{t}s_l),
\end{equation}
 for any natural number $t$.
\end{thm}

In this paper, by using a result of Babai, we give a lower bound for the maximum nullity of Cayley graph, $X_S(G)$, where $G=\langle a\rangle$ is a cyclic group, or $G=G_1\times \cdots\times G_t$ such that $G_1=\langle a\rangle$ is a cyclic group and $G_i$ is an arbitrary finite group, for some $2\leq i\leq t$, with determine the spectrum of Cayley graphs.
\section{Preliminaries}
For any positive integer $n$, define {\it M$\ddot{o}$bius} number, $\mu(n)$, as the sum of the primitive $n^{th}$ roots of unity. It has values in
$\{-1, 0, 1\}$ depending on the factorization of $n$ into prime factors.
\begin{enumerate}
  \item $\mu(1)=1$,
  \item $\mu(n)=0$, if $n$ has a squared factor,
  \item $\mu(n)=(-1)^{k}$,  if $n$ is a square-free with $k$ number of prime factors.
\end{enumerate}
Suppose that $k$ is a positive integer. The number of solutions of
$y_1+\cdots+y_r\equiv t\;\;(\text{mod}\;\; k)$,
where $y_1, \ldots, y_r$ and $t$ are belonged to the least non-negative residue system modulo $k$, is obtained in terms of the \text{von sterneck} function, $\Phi(n,k)$. In particular \text{von Sterneck} studied the case where the polynomial resulting from the expansion is reduced modulo a positive integer. This function is used in several equivalent forms; and in the form used by {\it $H\ddot{o}lder$} \cite{ho},
 \begin{equation}\label{von}
 \Phi(k,n)=\frac{\phi(n)}{\phi(n/(n,k))}\mu(n/(n,k)),
 \end{equation}
 where $k$ and $n$ are positive integers, $(n,k)$ is the greatest common divisor of $k$ and $n$, $\phi(n)$ is the {\it E$\ddot{u}$ler totient}, and $\mu (n)$ is the {\it M$\ddot{o}$bius}
 number. In the sequel, the following fundamental result is obtained by {\it $H\ddot{o}lder$}.
\begin{equation}\label{nic}
  \Phi(r,n)=\sum_{(r,n)=1} exp(2\pi irk/n).
\end{equation}
This properties was also studied by Von Sterneck in $1902$ \cite{von}, Nicol and Vandiver in $1954$ \cite{nic}, and Tom M. Apostol in $1972$ \cite{ap}.
\newline
Suppose that $B(k,n)=\big\{t\in \mathbb{N}\;:\; t\leq n\;,\; (t,n)=k\big\}$, and let $\omega=exp(2\pi i/n)$. Then the following function is called \text{Ramanujan} sum, and is denoted by $C(r,n)$.
\begin{equation}\label{ram}
\sum_{k\in B(1,n)}\omega^{kr},    0\leq r\leq n-1,
\end{equation}
 In \cite{Ramanujan}, it was obtained that $C(r,n)$ have only integral values, for some positive integers $r$ and $n$. Also, (\ref{nic}) and (\ref{ram}) state that $\Phi(r,n)=C(r,n)$.

\begin{lem}\label{r2} Suppose that $n>1$ and $d>1$ are two positive integers such that $d\mid n$. Also, let
 $B(d,n)=\{t\in \mathbb{N}\;:\; t\leq n\;,\; (t,n)=d\}$. Then
 \begin{enumerate}
   \item If $t\in B(d,n)$, then $C(t,n)=C(d,n)$,
   \item $|B(d,n)|=\phi(n/d)$.
 \end{enumerate}
\end{lem}
\noindent{\bf Proof.} The proof is straightforward.
 $\hfill\Box$
  \begin{thm}\label{maghsom}\cite{goldstein}
 For {\it E$\ddot{u}$ler totient} $\phi$ and positive integer $n$, we have $\sum\limits_{d|n}\phi(d)=n$.
\end{thm}
\begin{lem}\label{product}\cite{darafshe}
The irreducible character of $G\times H$ is $\chi \times \psi$ such
that $\chi$ and $\psi$ are the irreducible characters of G and H
respectively. the value of $\chi \times \psi$ for any $g\in G$ and
$h\in H$ is  $(\chi \times \psi)(g,h)=\chi(g)\psi(h).$
\end{lem}
\begin{lem}\label{character} \cite{darafshe} Let $G=\langle a\rangle$ be a cyclic group of order $n$. Then irreducible characters of $G$ are $\rho_j(a^{k})=\omega^{jk}$, where $j,k=0,1,\ldots,n-1$.
\end{lem}
\section{Main theorems}
 In the following theorem, we determine the spectrum of \text{Cayley} graph $X_S(G)$ whose $G$ is a cyclic group of order $n$. Here, we define $F(n_i)=(-1)^{n_i}\phi(n)/\phi(n_i)$, where $k_i$ is the number of prime factors in the decomposition of $n_i$.

\begin{thm}\label{cyclic1} Let $n$ be a positive integer and $D$ be its divisors set. Also, let $G=\langle a\rangle$ be a
cyclic group of order $n$ and $S=\{a^{i}\;:\;i\in B(1,n)\}$. Then
\begin{center}
  $spec(X_S(G))=\Bigg[\phi(n)^{1},0^{\sum\limits_{x\in X}\phi(x)},F(d_1)^{\phi(d_1)},...,F(d_t)^{\phi(d_t)}\Bigg]$,
\end{center}
where $X=\{d\in D\;:\;p^{2}\mid d\}$, for a prime $p$; and $d_i\in D\setminus X$, for some $1\leq i\leq t$.
\end{thm}
\noindent{\bf Proof.}
First, suppose that $n$ is a prime number. Thus  $X_S(G)$ is isomorphic to the complete graph $K_n$, and so
\begin{equation}
 spec(X_S(G))=\bigg[(p-1)^{1},-1^{(p-1)}\bigg].
\end{equation}
Now, consider the case in that $n$ is not prime. Let $\lambda_{\frac{n}{d_{i}}}$ be the eigenvalue of $X_S(G)$ corresponding to character of $\chi_{\frac{n}{d_{i}}}$, for some $d_i\in D$.
By Lemma \ref{Baba}, $\lambda_{\frac{n}{d_{i}}}=C(\frac{n}{d_{i}},n)$, and by the form used by {\it $H\ddot{o}lder$}  in (\ref{von}), we have
\begin{equation}
\lambda_{\frac{n}{d_{i}}}=\frac{\phi(n)}{\phi(d_i)}\mu(d_i).
\end{equation}
On the other hand, lemma \ref{r2} implies that the multiplicity of $\lambda_{\frac{n}{d_{i}}}$ is equal to $\phi(d_{i})$. If $d_i=1$, then $\lambda_n=\phi(n)$ with multiplicity $1$. Also, if $p^2\mid d_i$, then definition of {\it M$\ddot{o}$bius} number implies that $\lambda_{\frac{n}{d_{i}}}=0$. For other cases, $\lambda_{\frac{n}{d_{i}}}=F(d_i)$.
$\hfill\Box$
\newline

The following theorem which is proven by S. Akbari $et~al.$ \cite{yas}, help us to make a connection between the multiplicity of the eigenvalues of a graph $\mathcal{G}$ and its maximum nullity $M(\mathcal{G})$.
 \begin{thm}\label{max}\cite{yas}
  Let $\mathcal{G}$ be a graph of order $n$, and let $\lambda_i$ be its eigenvalue with respective multiplicity $n_i$. Then $M(\mathcal{G})\geq n_i$.
\end{thm}
As a result, Theorems \ref{cyclic1} and \ref{max}, state the following corollary.
\begin{cor}
Let $n$ be a positive integer and $D$ be its divisors set. Also, let $G=\big\langle a\big\rangle$ be a cyclic group of order $n$, and let $S=\{a^{i}\;:\;i\in B(1,n)\}$. For some prime $p$ and $d_i\in D$, the followings are established.
\begin{enumerate}
  \item If $n$ has a squared factor, then $M\big(X_S(G)\big)\geq max\bigg\{\sum\limits_{p^2\mid d_i}\phi(d_i), \phi(d_i)\bigg\}$.
  \item If $n$ is a square-free, then $M\big(X_S(G)\big)\geq \phi(d_i)$.
\end{enumerate}
\end{cor}
\begin{defn}
Let $G$ be a group, and let $S$ be a subset of $G$. Also, let $\Lambda=\big\{\chi_1, \ldots, \chi_k\big\}$ be the set of irreducible characters with degree $1$ of $G$. A character $\chi_i\in \Lambda$ is defined to be an $\ell$-index character of $G$, if has the same value $\ell$ on all letters in $S$; in other word, $\chi_i\in \Lambda$ is an $\ell$-index character of $G$ if $\chi(s_i)=\ell$, for all $s_i\in S$. In the sequel, An $\ell$-index number of $G$ is defined to be the number of $\ell$-index characters of $G$, and is denoted by $N_G(\ell)$.
\end{defn}
\begin{thm}\label{main}
Let $n$ be a positive integer whose divisors set is denoted by $D$. Also, let $G_1=\langle a\rangle$ be a cyclic group of order $n$, and let $S'=\big\{a^i\;:\; i\in B(1,n)\big\}$. Suppose that $G_2,\ldots, G_t$ are some arbitrary finite groups, and let $S_k$ is a subset of $G_k$, for some
$2\leq k\leq t$. If $S=\big\{(a^{i},\alpha_1, \ldots, \alpha_t)\;:\;a^i\in S', \alpha_k\in S_k\big\}$, then for some prime $p$ and $d_i\in D$, the followings are established.
\begin{enumerate}
  \item If $n$ has a square factor, then
   \begin{equation*}
   \begin{split}
      &M\big(X_S(G_1\times\cdots\times G_t)\big)\geq max  \\
        & \bigg\{\bigg(\prod_{i=2}^{t} \bigg
        (N_{G_i}\big(\ell_i\big)\big|S_i\big|\bigg)\bigg)\bigg(\sum\limits_{p^2\mid d_i}\phi\big(d_i\big)\bigg), \bigg(\prod_{i=2}^{t} \bigg
        (N_{G_i}\big(\ell_i\big)\big|S_i\big|\bigg)\bigg)\bigg(\phi\big(d_i\big)\bigg)\bigg\}.
   \end{split}
   \end{equation*}
  \item If $n$ is a square-free, then
  \begin{equation*}
    M\big(X_S(G_1\times\cdots\times G_t)\big)\geq \bigg(\prod\limits_{i=2}^{t} \bigg
        (N_{G_i}\big(\ell_i\big)\big|S_i\big|\bigg)\bigg)\bigg(\phi\big(d_i\big)\bigg).
  \end{equation*}
\end{enumerate}
\end{thm}
\noindent{\bf Proof.}
For some $2\leq k\leq t$, suppose that $\rho_{j_k}$ are the $\ell_k$-index irreducible characters with degree $1$ of $G_k$, and let $\chi_{\frac{n}{d_i}}$ be an irreducible character of $G_1$. Let $\lambda_{\frac{n}{d_i} j_1\ldots j_t}$ and $\lambda_{\frac{n}{d_i}}$ be the eigenvalues of $X_S(G_1\times\ldots\times G_t)$ and $X_{S'}(G_1)$ corresponding to irreducible characters of $\chi_{\frac{n}{d_i}}\times\rho_{j_1}\times\cdots\times \rho_{j_t}$ and $\chi_{\frac{n}{d_i}}$, respectively. Lemma \ref{product} implies that
  \begin{equation}
  \begin{split}
    \lambda_{\frac{n}{d_i} j_2\ldots j_t} & =\sum_{(g_1,\ldots,g_t)\in S}\bigg(\chi_{\frac{n}{d_i}}\times\rho_{j_2}\times\ldots\times\rho_{j_t}\bigg)(g_1,\ldots,g_t) \\
       & =\sum_{(g_1,\ldots,g_t)\in S}\bigg(\chi_{\frac{n}{d_i}}(g_1)\times\rho_{j_2}(g_2)\times\ldots\times\rho_{j_t}(g_t)\bigg)\\
       & =\bigg(\prod_{k=2}^{t} \bigg(N_{G_k}\big(\ell_k\big)\big|S_k\big|\bigg)\bigg)\sum_{s'\in S'}\bigg(\chi_{\frac{n}{d_i}}(s')\bigg).
  \end{split}
  \end{equation}
We have,
\begin{equation}
  \lambda_{\frac{n}{d_i} j_2\ldots j_t}=\bigg(\prod_{k=2}^{t} \bigg(N_{G_k}\big(\ell_k\big)\big|S_k\big|\bigg)\bigg)\big(\lambda_{\frac{n}{d_i}}\big).
\end{equation}
Hence, by Theorem \ref{cyclic1}, if $n$ is a square-free, then $\lambda_{\frac{n}{d_i} j_2\ldots j_t}$ is an eigenvalue of $X_S(G_1\times\cdots\times G_t)$ with multiplicity
\begin{equation}
  \bigg(\prod\limits_{k=2}^{t} \bigg(N_{G_k}\big(\ell_k\big)\big|S_k\big|\bigg)\bigg)\big(\phi(d_i)\big),
\end{equation}
 and if $n$ is divided by a prime number, then
$\lambda_{\frac{n}{d_i} j_2\ldots j_t}$ is an eigenvalue of $X_S(G_1\times\cdots\times G_t)$ with multiplicity
\begin{equation}
  \bigg(\prod\limits_{k=2}^{t} \bigg(N_{G_k}\big(\ell_k\big)\big|S_k\big|\bigg)\bigg)\big(\phi(d_i)\big),
\end{equation}
where $\lambda_{\frac{n}{d_i}}\neq 0$, or with multiplicity
\begin{equation}
  \bigg(\prod\limits_{k=2}^{t} \bigg(N_{G_k}\big(\ell_k\big)\big|S_k\big|\bigg)\bigg)\bigg(\sum\limits_{p^2\mid d_i}\phi(d_i)\bigg),
\end{equation}
where $\lambda_{\frac{n}{d_i}}\neq 0$. Therefore, theorem \ref{max}, completed the proof.
$\hfill\Box$
\newline

We apply the theorem to the dihedral groups $D_n$, with presentation
  \begin{equation}
  \big\langle a,b\;:\;a^n=b^2=e, (ab)^2=e\big\rangle.
  \end{equation}
For simplicity, we treat the case of odd $n$ only. For $n=2m+1$ there are $m$ irreducible character of degree $2$ and $2$ of degree $1$. Recall their character tables [\cite{darafshe}, chap.18], where $\omega$ denotes a primitive $n$th root of unity.
  \begin{table}[h!]
  \centering
\begin{tabular}{cccc}

&$a^i$&$a^ib$&\\ \hline
$\chi_j$ & $\omega^{jk}+\omega^{-jk}$ & $0$ & $(j=1,\ldots,m)$\\
$\chi_{m+1}$  & 1&  -1& \\
$\chi_{m+2}$  & 1&  1&

\end{tabular}
\caption{Character table of dihedral groups $D_n$, where $n=2m+1$.}\label{table1}
\end{table}
\newline
\begin{exam}
Let $G=\big\langle a\big\rangle$ be a group of order an odd $n$, and let $S=\big\{\big(a^i,a^j\big)\;:\; a^i\in G, a^j\in D_n, i,j\in B(1,n)\big\}$. Obviously,
$\chi_{m+1}$ and $\chi_{m+2}$ are two $1$-index irreducible characters of $D_n$, and so $N_{D_n}(1)=2$. Hence, by Theorem \ref{main}, we have
\begin{enumerate}
  \item If $n$ has a square factor, then
  \begin{equation*}
    M\big(Cay(G\times D_n:S)\big)\geq max\big\{ 2\phi(n)\big(\sum\limits_{p^2\mid d_i}\phi(d_i)\big), 2\phi(n)\big(\phi(d_i)\big)\big\}.
  \end{equation*}
  \item If $n$ is square-free, then $M\big(Cay(G\times D_n:S)\big)\geq 2\phi(n)\big(\phi(d_i)\big)$.
\end{enumerate}
\end{exam}

\noindent{\bf Acknowledgments.} This research  was  partially
supported by   Imam khomeini International University.


\begin{thebibliography}{99}

\bibitem{goldstein} W.W. Adams and L.J. Goldstein, {\em Introduction to Number Theory}, Prentice-Hall.inc, 1976.
 \bibitem{7} AIM Minimum Rank-Special Graphs Work Group, {\em Zero forcing sets and the minimum rank of graphs}, Linear Algebra and its Aplications, Vol. 428, 1628-1648, 2008.

\bibitem{yas} S. Akbari, E. Vatandoost and Y. Golkhandy Pour, {\em Maximum nullity and zero forcing number on cubic graphs}, arxiv:submit/1904157.


  \bibitem{ap} T.M. Apostol, {\em Arithmetical properties of generalized ramanujan sum}, pacific journal of mathematics, Vol. 41, N. 2, 281-293, 1972.

\bibitem{Babai} L. Babai, {\em Spectra of cayley graphs}, Journal of Combinatorial Theory, Vol. 27, 180-189, 1979.


 \bibitem{4} F. Barioli and S. Fallat, {\em On the minimum rank of the join of graphs and decomposable graphs}, Linear Algebra and its Applications, Vol. 421, 252-263, 2007.


\bibitem{2} A. Berman, S. Friedland, L. Hogben, U.G. Rothblum and B. Shader, {\em An upper bound for the minimum rank of a graph}, Linear
  Algebra and its Applications, Vol. 429, 1629-1638, 2008.


\bibitem{Biggs} N. Biggs, {\em Algebraic Graph Theory}, Cambridge University Press, London, 1993.


\bibitem{Doob} D.M. Cvetkovic, M.Doob and H. Sachs, {\em Spectra of graphs - Theory and applications}, 3rd edition, Johann Ambrosius Barth Verlag, Heidelberg-Leipzig, 1995.


\bibitem{von} R. Daublebsky Von Sterneck, {\em Ein Analogon zur additiven Zahlentheorie}, Sitzber Akad. Wiss. Wien, Math. Naturw. Klasse, vol. 111, 1567-1601, 1902.


  \bibitem{10} S.M. Fallat and L. Hogben, {\em The minimum rank of symmetric matrices described by a graph: A survey}, Linear Algebra and its Applications, Vol. 426, 558-582, 2007.


 \bibitem{1} L. Hogben, {\em Orthogonal representations, minimum rank, and graph complements}, Linear Algebra and its Applications, Vol. 428, 2560-2568, 2008.


       \bibitem{ho}  O. $H\ddot{o}lder$, {\em Zur Theorie der Kreisteilungsgleichung $K_m(x)=0$}, Prace Matematyczno Fizyczne, Vol. 43, 13-23, 1936.


\bibitem{darafshe} G. James and M. Liebeck, {\em Representations and Characters of groups}, Cambridge University Press, 1993.



\bibitem{Mc carty}  P.J. Maccarthy, {\em Introduction to arithmetical functions}, Universitext, Springer, New York, 1986.


       \bibitem{nic}   C.A. Nicol and H.S. Vandiver, {\em A Von Sterneck arithmetical function and restricted partitions with respect to a modulus}, Proc. Nat. Acad. Sci., Vol. 40, 825-835, 1954.



\bibitem{Ramanujan} S. Ramanujan, {\em On certain trigonometrical sums and their applications in the theory of numbers}, [Trans. Cambridge Philos. Soc. 22 (1918), no. 13, 259--276]. Collected papers of Srinivasa Ramanujan, 179-199, AMS Chelsea Publ., Providence, RI, 2000.








\end{thebibliography}
\end{document}